\documentclass[12pt]{amsart}

\newtheorem{theorem}{Theorem}[section]
\newtheorem{lemma}[theorem]{Lemma}

\newtheorem{corollary}[theorem]{Corollary}

\theoremstyle{definition}

\theoremstyle{remark}

\numberwithin{equation}{section}

\def\tr{\textup{tr}}

\def\pa{\partial}

\def\na{\nabla}

\newcommand\del{\delta}

\newcommand\ra{\longrightarrow}

\begin{document}

\title[N Square]{The Square of Nijenhuis Tensor and Its Vanishing  Results}

\author[Jun Ling]{Jun \underline{LING}}
\address{Department of Mathematics, Utah Valley University, Orem, UT 84058,USA\\
Visiting Fellow of Department of Mathematics, Princeton University, Princeton, NJ 08540, USA}
\email{lingju@uvu.edu\\ jl0026@princeton.edu}



\date{October 18, 2019}



\begin{abstract}
	We give the strong form and the weak form of the square of Nijenhuis tensor,
	and some vanish results of the square.
\end{abstract}

\maketitle

\section{Introduction}
 A complex manifold is a (real) smooth manifold, since the complex (holomorphic) structure of the complex manifold yields a smooth differentiable structure for the underlying real manifold. The converse is not true in general. 

\vspace{0.1in}

Whether or not a smooth manifold,
say spheres, is a complex manifold as well becomes an important issue. More specifically,
given a complex manifold, its complex (holomorphic) structure gives a smooth differentiable structure, and in addition, gives an underlying almost-complex structure,  an endomorphism $J$ on the tangent bundle of underlying smooth (real) manifold with $J^2=-1$. Conversely, can a smooth manifold $M$ with an almost-complex structure $J$ be a complex manifold such  that $J$ is the underlying almost-complex structure of the complex structure of the complex manifold? This is an interesting question and it attracts many mathematicians' attention.   

\vspace{0.1in}

There have been progress and very interesting results. It is known that a  smooth manifold with complex structure must be even dimensional and orientable. 

\vspace{0.1in}

In sphere cases, it is known that $S^2$ is a complex manifold.
For high dimensional spheres other than 6-sphere, there have been a lot results of non-existences  of complex structure, 
ch. Ehresmann \cite{Ehr}, H. Hopf \cite{Hopf},
A. Kirchhoff \cite{Kir},
Eckmann-Frohlicher \cite{EF}, Ehresmann-Liberman \cite{EL}, Borel and Serre \cite{B-S} 
The above is only very incomplete short list of results on this problem.
Hirzebruch \cite{Hir} in 1954 and Liberman\cite{Lib} in 1955, Yau \cite{Yau} asked whether or not  there is a complex structure on 6-sphere? 

\vspace{0.1in}

When we study the problem whether  a given smooth manifold is a complex manifold or not, one 
 approach is to begin with an almost-complex structure and to check further if the almost-complex structure can be integrated in to the complex structure of a complex manifold. When checking, one applies the  criterion of   Newlander-Nirenberg\cite{NN}, which says that an almost-complex structure
$J$ can be  integrated into a complex structure if and only if the Nijenhuis tensor $N$ of the almost-complex structure $J$ vanishes. 
Therefore the study of Nijenhuis tensor becomes interesting and crucial here. 
 
\vspace{0.1in}

In this paper, we write down some calculations on the "square" of Nijenhuis tensor, which is  in weak and strong  meaning,  and vanishing results of the the square of Nijenhuis tensor, in weak and strong  meaning. The motivation for us to look at the square of Nijenhuis tensor and vanishing results of square for certain geometric objects is we try to classify the image of Nijenhuis tensor, and wish to ralate those classes to cohomology and seek properties and relations for almost-complex structure and complex structures.

\vspace{0.1in}

We state the settings, definitions and  main results in Section \ref{sec-main-results}, In Section \ref{sec-proof-T-is0} we prove  that  the weak form of square of Nijenhuis tensor is zero for general almost-complex manifolds. That section is the revision of the author's another work about one and one half year ago. In Section \ref{sec-property-L-ell} we  study properties of the square of the Nijenhuis tensor and the recapture of the strong form of square from the intermediate form of square. In Section \ref{sec-proof-other-main results} we prove the  other main result.

\vspace{0.2in}

\section{Settings and Main Results}\label{sec-main-results}
	Let $J$ be an almost-complex structure on a smooth $n$-manifold $M$, that is smooth $1-1$ tensor field with $J^2=-1$. $N=N_J$ is the Nijenhuis tensor of the almost-complex structure $J$, given by the following equation
\begin{equation}\label{Nijen}
N(X,Y)=[JX,JY]-J[X,JY]-J[JX,Y]-[X,Y]
\end{equation}
for all smooth vector fields $X$ and $Y$ on manifold $M$, where $[X,Y]=XY-YX$ is the Lie bracket. It is easy to see that $N$ is a tensor, so the value of $N(X,Y)(p)$ depends on vectors  $X(p)$ and $Y(p)$ (for $p\in M$) only, not their extensions.

We consider the \textbf{strong form} $N^2$ of the square of Nijenhuis  tensor $N$:
\[
N^2:\Gamma(TM)\times \Gamma(TM)\times \Gamma(TM)\ra \Gamma(TM),
\]
\begin{equation}\label{N-sq-strong}
N^2(X,Y,Z):=N\Big(N(X,Y),Z\Big).
\end{equation}
Note that $N^2$ depends on $N$ only.

\vspace{0.1in}

On the smooth manifold $M$, we take a Riemannian metric $g$. Using metric $g$ we consider the following tensor $L=L_g$. For all smooth vector fields $X$ and $Y$ on manifold $M$, we define $
L_g:=L$ by
\begin{equation}\label{L-def}
L(X,Z,Y,W):=\frac14\Big\{\left\langle JN\Big(N(X,Z),Y\Big),W\right\rangle_g
\end{equation}
\[
\quad\qquad\qquad\qquad\qquad+\left\langle JN\Big(N(Y,Z),X\Big),W\right\rangle_g
\]
\[
\quad\qquad\qquad\qquad\qquad+\left\langle JN\Big(N(X,W),Y\Big),Z\right\rangle_g
\]
\[
\qquad\qquad\qquad\qquad\qquad+\left\langle JN\Big(N(Y,W),X\Big),Z\right\rangle_g\Big\},
\]
here and below $\langle\cdot,\cdot\cdot \rangle_g =g(\cdot,\cdot\cdot)$. 

\vspace{0.1in}

We define the
\textbf{intermediate form} $\ell=\ell_g$ of the square of the Nijenhuis tensor.
\begin{equation}\label{ell-def}
\ell(X,Z):=L(X,Z,X,Z).
\end{equation}
Noticing that by the symmetry in  (\ref{L-def}), we actually have
\[
\ell(X,Z)=\left\langle JN\Big(N(X,Z),X\Big),Z\right\rangle_g.
\]

\vspace{0.1in}

Consider function $T$ on $M$ obtained by the following equation, 
\[
T:=\tr_{1,3}\Big(\tr_{2,4}\big(L(\textup{arg}_1,\textup{arg}_2, \textup{arg}_3, \textup{arg}_4\big)\Big),
\]
where $\tr_{1,3}$ is the trace of the 1st argument $\textup{arg}_1,$ and 3rd argument $\textup{arg}_3$.  $\tr_{2,4}$ has similar meaning.

\vspace{0.1in}

We call $T$ the \textbf{weak form} of $N^2$.

\vspace{0.1in}

Note that though both $L$ and $\ell$ depend on metric $g$, $T$ does not. The following
Theorem \ref{T-is0-thm} shows $T$ is independent of metric $g$.

\begin{theorem}\label{T-is0-thm}
	For smooth manifold $M$ with almost-complex structure $J$,
the weak form $T$ of $N^2$ is 
independent of Riemannian metric $g$ and is zero function:
\[
T\equiv0.
\]
\end{theorem}

\vspace{0.1in}

\begin{theorem}\label{NN-is0-thm}
	If for some Riemannian metric $g$,  the intermediate form $\ell_g\equiv 0$, then the strong form  $N^2$ is zero:
	\begin{equation}\label{NNis0}
	N^2\equiv 0.
	\end{equation}
	Again the result is independent of Riemannian metric $g$.
\end{theorem}

\vspace{0.2in}

\section{Vanishing of the Weak Form $T$ of $N^2$} \label{sec-proof-T-is0}
In this section we present a proof for Theorem \ref{T-is0-thm}.

\proof
Take a local normal coordinates  $\{x^i\}_{i=1}^n$ at $x\in M$ with $\{e_i=\dfrac{\pa}{\pa x^i}\}_{i=1}^n$ , let $g_{ij}=g(e_i,e_j)$, as usual.
Take the trace of the  first argument and the third argument and then take the trace
of the second and the fourth argument, then we have
\begin{equation}\label{T-def}
T=\sum_{i,j,k,l}g^{ij}g^{kl}L(e_i,e_k,e_j,e_l),
\end{equation}
where  matrix $(g^{ij})=(g_{ij})^{-1}$. 
We write $\pa_i=\frac{\pa}{\pa x^i}$, $\na_i=\na_{\frac{\pa}{\pa x^i}}$
$\langle\cdot,\cdot\cdot \rangle_g =g(\cdot,\cdot\cdot)$ for convenience.

\vspace{0.1in}

To prove Theorem \ref{T-is0-thm} is to show that 
\begin{equation}\label{T-is0-eq}
T=0.
\end{equation}

\vspace{0.1in}

Recall from (\ref{L-def})
\[
L(X,Z,Y,W):=\frac14\Big\{\left\langle JN\Big(N(X,Z),Y\Big),W\right\rangle_g
\]
\[
\quad\qquad\qquad\qquad\qquad+\left\langle JN\Big(N(Y,Z),X\Big),W\right\rangle_g
\]
\[
\quad\qquad\qquad\qquad\qquad+\left\langle JN\Big(N(X,W),Y\Big),Z\right\rangle_g
\]
\[
\qquad\qquad\qquad\qquad\qquad+\left\langle JN\Big(N(Y,W),X\Big),Z\right\rangle_g\Big\},
\]
where $N=N_J$ is the one in (\ref{Nijen}). We calculate at a point $x$ with local normal coordinates  $\{x^i\}_{i=1}^n$ with $\{e_i=\dfrac{\pa}{\pa x^i}\}_{i=1}^n$. So at the point $x$,  $g_{ij}=\del_{ij}$,
$\del_{ii}=1$, $\del_{ij}=0$ if $i\not=j$. The following computations are at point $x$.

\vspace{0.1in}

Note that $L(X,Z,X,Z)=\left\langle JN\Big(N(X,Z),X\Big),Z\right\rangle_g$.
Therefore the trace is
\[
T=\sum_{i,k}\ell(e_i,e_k)=\sum_{i,k}L(e_i,e_k,e_i,e_k)=\sum_{i,k}\left\langle JN\Big(N(e_i,e_k),e_i\Big),e_k\right\rangle_g
=\sum_{i,k,r,s}N_{ik}^rN_{ri}^sJ_s^k,
\]
where $N_{ij}^k$ are given by $N(\pa_i,\pa_j)=N_{ij}^k\pa_k$. 
It is easy to see
\begin{equation}\label{N-ijk}
N_{ij}^k=J_i^k(\pa_pJ^k_j-\pa_jJ^k_p)-J_j^k(\pa_pJ^k_i-\pa_iJ^k_p).
\end{equation}
We write 
\[J_i:=\sum_{i=1}^nJ_i^j\pa_j\] 
for convenience. In the following we omit sum symbol and \textbf{take sum over double indices}.

\vspace{0.1in}

By (\ref{N-ijk}), we have 
\begin{equation}\label{T-eq}
T=N_{ik}^rN_{ri}^sJ_s^k=J_s^kN_{ik}^rN_{ri}^s
\end{equation}
\[
=J_s^k\{J_i^p(\pa_pJ^r_k-\pa_kJ^r_p)-J_k^p(\pa_pJ^r_i-\pa_iJ^r_p)\}
\{J_r^q(\pa_qJ^s_i-\pa_iJ^s_q)-J_i^q(\pa_qJ^s_r-\pa_rJ^s_q)\}
\]
\[
=\{J_s^kJ_i^p(\pa_pJ^r_k-\pa_kJ^r_p)-J_s^kJ_k^p(\pa_pJ^r_i-\pa_iJ^r_p)\}
\{J_r^q(\pa_qJ^s_i-\pa_iJ^s_q)-J_i^q(\pa_qJ^s_r-\pa_rJ^s_q)\}
\]
\[
=\{J_s^kJ_i^p\pa_pJ^r_k-J_s^kJ_i^p\pa_kJ^r_p
+\del_s^p(\pa_pJ^r_i-\pa_iJ^r_p)\}
\{J_r^q(\pa_qJ^s_i-\pa_iJ^s_q)-J_i^q(\pa_qJ^s_r-\pa_rJ^s_q)\}
\]
\[
=\{J_s^kJ_i^p\pa_pJ^r_k-J_s^kJ_i^p\pa_kJ^r_p
+\pa_sJ^r_i-\pa_iJ^r_s\}
\{J_r^q\pa_qJ^s_i-J_r^q\pa_iJ^s_q-J_i^q\pa_qJ^s_r+J_i^q\pa_rJ^s_q\}
\]
\[
=\{J_s^k\cdot J_iJ^r_k-J_i^p\cdot J_sJ^r_p
+\pa_sJ^r_i-\pa_iJ^r_s\}
\{J_rJ^s_i-J_r^q\cdot \pa_iJ^s_q-J_iJ^s_r+J_i^q\cdot \pa_rJ^s_q\}
\]

\[
\stackrel{\textup{def}}{=}\{a+b+c+d\}\{(1)+(2)+(3)+(4)\}
\]
\[
=a(1)+a(2)+a(3)+a(4)
\]
\[
+b(1)+b(2)+b(3)+b(4)
\]
\[
+c(1)+c(2)+c(3)+c(4)
\]
\[
+d(1)+d(2)+d(3)+d(4)
\]
\[
:=J_s^k\cdot J_iJ^r_k\cdot J_rJ^s_i-J_r^qJ_s^k\cdot J_iJ^r_k\cdot \pa_iJ^s_q-J_s^k\cdot J_iJ^r_k\cdot J_iJ^s_r+J_i^qJ_s^k\cdot J_iJ^r_k\cdot \pa_rJ^s_q
\]
\[
-J_i^p\cdot J_sJ^r_p\cdot J_rJ^s_i+J_r^qJ_i^p\cdot J_sJ^r_p\cdot \pa_iJ^s_q+J_i^p\cdot J_sJ^r_p\cdot J_iJ^s_r-J_i^qJ_i^p\cdot J_sJ^r_p\cdot \pa_rJ^s_q
\]
\[
+J_rJ^s_i\cdot\pa_sJ^r_i-J_r^q\cdot \pa_iJ^s_q\cdot\pa_sJ^r_i-J_iJ^s_r\cdot\pa_sJ^r_i+J_i^q\cdot \pa_rJ^s_q\cdot\pa_sJ^r_i
\]
\[
-J_rJ^s_i\cdot\pa_iJ^r_s+J_r^q\cdot\pa_iJ^r_s\cdot \pa_iJ^s_q+J_iJ^s_r\cdot\pa_iJ^r_s-J_i^q\cdot\pa_rJ^s_q\cdot\pa_iJ^r_s
\]
where
in the second term of the first factor we have used 
\[
J_s^kJ_k^p=-\del_s^p.
\]

\vspace{0.1in}

In the following, we  will use this relation exclusively without further explanation. We compute terms in (\ref{T-eq}).

\[a(1)=J_s^k\cdot J_iJ^r_k\cdot J_rJ^s_i=J_l^j\cdot J_iJ^k_j\cdot J_kJ^l_i,
\]
and
\[b(3)
=J_i^p\cdot J_sJ^r_p\cdot J_iJ^s_r
\]
\[
=J_i^p\cdot J_iJ^s_r\cdot  J_sJ^r_p
\]
\[
=J_i^l\cdot J_iJ^k_j \cdot J_kJ^j_l
\]
\[
=J_i^l\cdot J_iJ^k_j \cdot J_k^p\cdot\pa_pJ^j_l
\]
\[
=J_k^p\cdot J_iJ^k_j \cdot J_i^l\cdot\pa_pJ^j_l
\]
\[
=-J_k^p\cdot J_iJ^k_j \cdot J^j_l \cdot\pa_p J_i^l
\]
\[
=-J_k^pJ^j_l\cdot J_iJ^k_j \cdot \pa_p J_i^l
\]
\[
=-J^j_l\cdot J_iJ^k_j \cdot J_k^p\cdot\pa_p J_i^l
\]
\[
=-J^j_l\cdot J_iJ^k_j \cdot J_k J_i^l.
\]
Therefore we have
\[
a(1)+b(3)=0.
\]

\vspace{0.1in}

\[
a(2)=-J_r^qJ_s^k\cdot J_iJ^r_k\cdot \pa_iJ^s_q
\]
\[
=J^s_qJ_s^k\cdot J_iJ^r_k\cdot \pa_iJ_r^q
\]
\[=-J_iJ^r_q\cdot \pa_iJ_r^q
\]
\[=- J_iJ^k_j\cdot \pa_iJ_k^j,
\]
and
\[
d(3)=J_iJ^s_r\cdot\pa_iJ^r_s
\]
\[
=J_iJ^k_j\cdot\pa_iJ^j_k.
\]
Therefore we have
\[
a(2)+d(3)=0.
\]

\vspace{0.1in}

\[
a(3)=-J_l^j\cdot J_iJ^k_j\cdot J_iJ^l_k
\]
\[
=-J_t^j\cdot J_iJ^k_j\cdot J_iJ^t_k
\]
\[
=-J_t^j\cdot J_iJ^l_j\cdot J_iJ^t_l
\]
\[
=-J_t^k\cdot J_iJ^l_k\cdot J_iJ^t_l
\]
So
\[
a(3)=-J_t^k\cdot J_iJ^l_k\cdot J_iJ^t_l
\]
\[
=-J_t^k\cdot J_pJ^l_k\cdot J_pJ^t_l
\]
\[
=-J_t^k\cdot J_p^i\pa_iJ^l_k\cdot J_p^j\pa_jJ^t_l
\]
\[
=-\pa_i(J_t^kJ_p^iJ_p^j\cdot J^l_k\cdot \pa_jJ^t_l)
\]
\[
+\pa_iJ_t^k\cdot J_p^iJ_p^j\cdot J^l_k\cdot \pa_jJ^t_l
\]
\[
+J_t^k\pa_i J_p^i\cdot J_p^jJ^l_k\cdot \pa_jJ^t_l
\]
\[
+J_t^kJ_p^i\pa_iJ_p^j\cdot J^l_k\cdot \pa_jJ^t_l
\]
\[
+J_t^kJ_p^iJ_p^j\cdot J^l_k\cdot \pa_j\pa_iJ^t_l
\]
\[
=+\pa_i(J_p^iJ_p^j\cdot \pa_jJ^l_l)
\]
\[
+J_p^iJ_p^j\cdot\pa_iJ_t^k\cdot J^l_k\cdot \pa_jJ^t_l
\]
\[
-J_p^j\cdot\pa_i J_p^i\cdot \pa_jJ^l_l
\]
\[
-J_p^i\cdot \pa_iJ_p^j\cdot  \pa_jJ^l_l
\]
\[
-J_p^iJ_p^j\cdot \pa_j\pa_iJ^l_l
\]
\[
=J^l_kJ_p^iJ_p^j\cdot\pa_iJ_t^k\cdot\pa_jJ^t_l
\]
\[
=J^l_kJ_p^iJ_p^j\cdot\pa_jJ^t_l\cdot\pa_iJ_t^k
\]
\[
=J^b_aJ_p^iJ_p^j\cdot\pa_jJ^c_b\cdot\pa_iJ_c^a
\]
\[
=J^k_tJ_p^iJ_p^j\cdot\pa_jJ^l_k\cdot\pa_iJ_l^t
\]
\[
=J^k_tJ_p^uJ_p^v\cdot\pa_vJ^l_k\cdot\pa_uJ_l^t
\]
\[
=J^k_tJ_p^vJ_p^u\cdot\pa_vJ^l_k\cdot\pa_uJ_l^t
\]
\[
=J^k_tJ_p^iJ_p^j\cdot\pa_iJ^l_k\cdot\pa_jJ_l^t
\]
\[
=-a(3),
\]
that is,
\[
a(3)=-a(3).
\]
Therefore we have
\[
a(3)=0.
\]

\vspace{0.1in}

\[
a(4)=J_i^qJ_s^k\cdot J_iJ^r_k\cdot \pa_rJ^s_q
\]
\[
=-J^s_qJ_s^k\cdot J_iJ^r_k\cdot \pa_rJ_i^q
\]
\[
=J_iJ^r_k\cdot \pa_rJ_i^k
\]
\[
=J_iJ^k_j\cdot \pa_kJ_i^j,
\]
where
\[
J_i^q \cdot\pa_rJ^s_q=-\pa_rJ_i^q \cdot J^s_q,
\]
and
\[J^s_qJ_s^k=-\del^s_k,
\]

\[
c(3)=-J_iJ^s_r\cdot\pa_sJ^r_i=-J_iJ^k_j\cdot\pa_kJ^j_i.
\]
Therefore we have
\[
a(4)+c(3)=0.
\]

\vspace{0.1in}

\[
d(1)=-J_rJ^s_i\cdot\pa_iJ^r_s
\]
\[
=-J_jJ^k_i\cdot\pa_iJ^j_k,
\]
where we changed index from $r$ ro $j$.
\[
c(2)=-J_r^q\pa_iJ^s_q\cdot\pa_sJ^r_i
\]
\[
=-\pa_sJ^r_i\cdot \pa_iJ^s_q\cdot J_r^q
\]
\[
=J^s_q\cdot \pa_sJ^r_i\cdot \pa_iJ_r^q
\]
\[
=J_qJ^r_i\cdot \pa_iJ_r^q
\]
\[
=J_jJ^k_i\cdot \pa_iJ_k^j,
\]
where we used
\[
\pa_iJ^s_q\cdot J_r^q=\pa_i(J^s_q\cdot J_r^q)-J^s_q\cdot \pa_iJ_r^q
\]
\[=\pa_i(\del^s_r)-J^s_q\cdot \pa_iJ_r^q
\]
\[=-J^s_q\cdot \pa_iJ_r^q.
\]
Therefore
\[
d(1)+c(2)=0.
\]

\vspace{0.1in}

\[
b(4)=-J_i^qJ_i^p\cdot J_sJ^r_p\cdot \pa_rJ^s_q
\]
\[
=-J_i^p\cdot J_sJ^r_p\cdot J_i^q\cdot\pa_rJ^s_q
\]
\[
=J_i^p\cdot J_sJ^r_p\cdot J^s_q\cdot\pa_rJ_i^q
\]
\[
=J_i^p\cdot J_s^t\cdot\pa_tJ^r_p\cdot J^s_q\cdot\pa_rJ_i^q
\]
\[
=J_i^p\cdot J_s^t\cdot J^s_q\cdot\pa_tJ^r_p\cdot\pa_rJ_i^q
\]
\[
=-J_i^p\cdot\del_q^t\cdot\pa_tJ^r_p\cdot\pa_rJ_i^q
\]
\[
=-J_i^p\cdot\pa_qJ^r_p\cdot\pa_rJ_i^q
\]
\[
=J^r_p\cdot\pa_qJ_i^p\cdot\pa_rJ_i^q
\]
\[
=J^r_p\cdot\pa_rJ_i^q\cdot \pa_qJ_i^p
\]
\[
=J^r_j\cdot\pa_rJ_i^k\cdot \pa_kJ_i^j
\]
\[
=J_jJ_i^k\cdot \pa_kJ_i^j
\]
\[=c(1),
\]
that is,
\[
b(4)=c(1).
\]

\vspace{0.1in}

\[b(1)=-J_i^p\cdot J_sJ^r_p\cdot J_rJ^s_i
\]
\[
=-J_i^p\cdot J_rJ^s_i\cdot J_sJ^r_p
\]
\[
=-J_i^l\cdot J_jJ^s_i\cdot J_sJ^j_l\]
\[
=-J_i^l\cdot J_jJ^k_i \cdot 
J_kJ^j_l
\]
\[
=-J_i^l\cdot J_jJ^k_i \cdot 
J_k^p\cdot\pa_pJ^j_l
\]
\[
=-J_k^p\cdot J_jJ^k_i \cdot 
J_i^l\cdot\pa_pJ^j_l
\]
\[
=J_k^p\cdot J_jJ^k_i \cdot J^j_l
\cdot\pa_pJ_i^l
\]
\[
=J^j_l\cdot J_jJ^k_i \cdot 
J_k^p\pa_pJ_i^l
\]
\[
=J^j_l\cdot J_jJ^k_i \cdot J_kJ_i^l
\]
\[
=J^j_l\cdot J_j^p\pa_pJ^k_i \cdot J_kJ_i^l
\]
\[
=-\del_l^p\pa_pJ^k_i \cdot J_kJ_i^l
\]
\[
=-\pa_lJ^k_i \cdot J_kJ_i^l
\]
\[
=- J_kJ_i^l\cdot \pa_lJ^k_i
\]
\[
=- J_jJ_i^k\cdot \pa_kJ^j_i.
\]

\[
c(1)=J_rJ^s_i\cdot\pa_sJ^r_i
\]
\[
=J_jJ^k_i\cdot\pa_kJ^j_i
\]
Therefore
\[
b(1)+b(4)=(b(1)+c(1))=0.
\]

\vspace{0.1in}

\[
b(2)=J_r^qJ_i^p\cdot J_sJ^r_p\cdot \pa_iJ^s_q
\]
\[
=-J_i^pJ^s_q\cdot J_sJ^r_p\cdot \pa_iJ_r^q
\]
\[
=-J_i^pJ^s_q\cdot J_s^k\pa_kJ^r_p\cdot \pa_iJ_r^q
\]
\[
=J_i^p\pa_kJ^r_p\cdot \pa_iJ_r^k
\]
\[
=J_i^p\pa_jJ^r_p\cdot \pa_iJ_r^j
\]
\[
=-J^l_k\pa_jJ_i^k\cdot \pa_iJ_l^j
\]
\[
=J_l^j\pa_jJ_i^k\cdot \pa_iJ^l_k
\]
\[
=J_lJ_i^k\cdot \pa_iJ^l_k
\]
\[
=J_jJ_i^k\cdot \pa_iJ^j_k
\]
and
\[
d(4)=-J_i^q\pa_rJ^s_q\cdot\pa_iJ^r_s
\]
\[
=J^s_q\cdot \pa_rJ_i^q\cdot\pa_iJ^r_s
\]
\[
=J^s_q\cdot \pa_jJ_i^q\cdot\pa_iJ^j_s
\]
\[
=-J^j_s\cdot \pa_jJ_i^q\cdot  \pa_i J^s_q
\]
\[
=-J_sJ_i^q\cdot  \pa_i J^s_q=-J_jJ_i^k\cdot  \pa_i J^j_k,
\]
Therefore we have
\[
b(2)+d(4)=0.
\]

\vspace{0.1in}

\[
c(1)=J_rJ^s_i\cdot\pa_sJ_i^r,
\]

\[
c(4)=J_i^q\pa_rJ^s_q\cdot\pa_sJ^r_i
\]
\[
=-J^r_i\pa_rJ^s_q\cdot\pa_sJ_i^q
\]
\[
=-J_iJ^s_q\cdot\pa_sJ_i^q
\]
\[
=-J_iJ^k_j\cdot\pa_kJ_i^j.
\]
and
\[
c(4)=J_i^q\pa_rJ^s_q\cdot\pa_sJ^r_i
\]
\[
=-J^s_q\pa_rJ_i^q\cdot\pa_sJ^r_i
\]
\[
=-J^s_q\pa_sJ^r_i\cdot\pa_rJ_i^q
\]
\[
=-J_qJ^r_i\cdot\pa_rJ_i^q
\]
\[
=-J_qJ^s_i\cdot\pa_sJ_i^q
\]
\[
=-J_rJ^s_i\cdot\pa_sJ_i^r
\]
\[
=-c(1).
\]
Therefore we have
\[
c(1)+c(4)=0.
\]

\vspace{0.1in}
Taking all the above  into (\ref{T-eq}), we have

\[
T=N_{ik}^rN_{ri}^sJ_s^k
\]
\[
=-J_s^k\cdot J_iJ^r_k\cdot J_iJ^s_r+J_s^k\cdot J_iJ^r_k\cdot J_rJ^s_i+J_i^p\cdot J_sJ^r_p\cdot J_iJ^s_r-J_i^p\cdot J_sJ^r_p\cdot J_rJ^s_i
\]
\[
-J_r^qJ_s^k\cdot J_iJ^r_k\cdot \pa_iJ^s_q+J_i^qJ_s^k\cdot J_iJ^r_k\cdot \pa_rJ^s_q+II3
+J_iJ^s_r\cdot\pa_iJ^r_s-J^l_k\pa_jJ_i^k\cdot \pa_iJ_l^j
\]
\[
-J_i^qJ_i^p\cdot J_sJ^r_p\cdot \pa_rJ^s_q-J_iJ^s_r\cdot\pa_sJ^r_i+J_rJ^s_i\cdot\pa_sJ^r_i
\]
\[
+J_r^q\pa_iJ^s_q\cdot\pa_iJ^r_s-J_jJ_i^k\cdot  \pa_i J^j_k
+IV3+J_i^q\pa_rJ^s_q\cdot\pa_sJ^r_i.
\]
\[
=a(3)+a(1)+b(3)+b(1)
\]
\[
+a(2)+a(4)+d(1)+d(3)+b(2)
\]
\[
+b(4)+c(3)+c(1)
\]
\[
+d(2)+d(4)+c(2)+c(4)
\]
\[
=d(2)
\]
\[
=0,
\]
since 
\[
d(2)=0.
\]

\vspace{0.1in}

In fact,
\[
d(2)
=J_l^j\pa_iJ^k_j\cdot\pa_iJ^l_k
\]
\[
=\pa_i(J_l^jJ^k_j\cdot\pa_iJ^l_k)
\]
\[
-\pa_iJ_l^j\cdot J^k_j\cdot\pa_iJ^l_k
\]
\[
-J_l^jJ^k_j\cdot\pa_i\pa_iJ^l_k
\]
\[
=-\pa_i(\pa_iJ^k_k)
\]
\[
-\pa_iJ_l^j\cdot J^k_j\cdot\pa_iJ^l_k
\]
\[
+\pa_i\pa_iJ^k_k
\]
\[
=-\pa_iJ_l^j\cdot J^k_j\cdot\pa_iJ^l_k
\]
\[
=-J^k_j\cdot \pa_iJ^l_k\cdot\pa_iJ_l^j
\]
Therefore
\[
d(2)=-d(2),
\]
\[
d(2)=0.
\]

\vspace{0.1in}

Therefore
\[
T:=N_{ik}^rN_{ri}^sJ_s^k=0.
\]
\qed
 
 \vspace{0.2in}
 
\section{Properties of $L$ and $\ell$, Recapture of $L$ From $\ell$ }\label{sec-property-L-ell}

We study some properties $L$ and $\ell$, which are needed in the our major results, though themselves make independent sense.
Tensor $L$ have the following properties.
\begin{lemma}[Symmetry and Asymmetry]\label{L-is-0-lem}
	\[
	L(X,Z,Y,W)=	L(Y,Z,X,W),
	\]
	\[
	L(X,Z,Y,W)=	L(X,W,Y,Z),
	\]
	\[	
	L(X,Z,X,Z)=\left\langle JN\Big(N(X,Z),X\Big),Z\right\rangle_g
	\]
	\[
	L(X,Z,X,Z)
	=-L(Z,X,X,Z).
	\]
	\[
	L(X,Z,X,Z)
	=L(Z,X,Z,X).
	\]
	for all smooth vector fields $X, Y, Z, W$ on manifold $M$.
\end{lemma}

\proof
All are from the definition (\ref{L-def}).
The last equality is Corollary \ref{ell-symm-coro}.
\qed

\vspace{0.1in}

\begin{lemma}
	We have he following results.
	\[
	N\Big(N(JX,Z),Y\Big)=JN\Big(N(X,Z),Y\Big)
	\]
	and
	\[
	N\Big(N(JX,Z),JY\Big)=N\Big(N(X,Z),Y\Big).
	\]
\end{lemma}
\proof

Direct computations shows 
\[
N(X,JY)
\]
\[
=[JX,JJY]-J[X,JJY]-J[JX,JY]-[X,JY]
\]
\[
=[JX,-Y]-J[X,-Y]-J[JX,JY]-[X,JY]
\]
\[
=-[JX,Y]+J[X,Y]-J[JX,JY]-[X,JY]
\]
\[
=-J[JX,JY]-[X,JY]-[JX,Y]+J[X,Y].
\]
\[
=-J\big([JX,JY]-J[X,JY]-J[JX,Y]-[X,Y]\big)
\]
\[
=-JN(X,Y),
\]
that is
\[
N(X,JY)=-JN(X,Y).
\]
Thus
\[
N(JX,Y)=-N(Y,JX)=JN(Y,X)=-JN(X,Y),
\]

\[
N(JX,JY)=-JN(JX,Y)=(-J)(-J)N(X,Y)=-N(X,Y),
\]
and
\[
N(X,JX)=-JN(X,X)=0.
\]
Therefore by he results above, we have
\[
N\Big(N(JX,Z),Y\Big)=
N\Big(-JN(X,Z),Y\Big)=JN\Big(N(X,Z),Y\Big)
\]
and
\[
N\Big(N(JX,Z),JY\Big)=
N\Big(-JN(X,Z),JY\Big)=N\Big(N(X,Z),Y\Big).
\]

\qed

\vspace{0.1in}

\begin{theorem}[Dependence of the Plane Only]\label{ellXZ-transition-thm}
	Suppose that $\{X,Z\}$ and $\{X',Z'\}$ are two pairs of linearly independent sets of vector fields spaning the same vector space at every point in an open set $U$. Then
	\[
	\ell(X',Z')=\left(\frac{\pa\{X',Z'\}}{\pa\{X,Z\}}\right)^2\ell(X,Z), \quad\textup{on }U
	\]
	and
	\[
	\frac{\ell(X',Z')}{\langle X',X'\rangle_g\langle Z',Z'\rangle_g-(\langle X',Z'\rangle_g)^2}
	=\frac{\ell(X,Z)}{\langle X,X\rangle_g\langle Z,Z\rangle_g-(\langle X,Z\rangle_g)^2}, \quad\textup{on }U,
	\]
	where $\frac{\pa\{X',Z'\}}{\pa\{X,Z\}}$ is the Jacobian.
\end{theorem}

\vspace{0.1in}

\begin{corollary}\label{ell-symm-coro}
	\[
	\ell(X,Z)=\ell(Z,X).
	\]
\end{corollary}
\proof[Proof of Corollary \ref{ell-symm-coro}] By Theorem \ref{ellXZ-transition-thm}
\[
\ell(Z,X)=(-1)^2\ell(X,Z)=\ell(X,Z).
\]
\qed

\vspace{0.1in}

\proof[Proof of Theorem \ref{ellXZ-transition-thm}] Calculate at one point $p\in U$.
If $X'=aX+bZ$ and $Z'=cX+dZ$ with constants $a,b,c,d$ and $ad-bc\not=0$, then
\[
\ell(X',Z')=
\left\langle JN\Big(N(X', Z'),X'\Big),Z'\right\rangle_g
\]
\[
=
\left\langle JN\Big(N(aX+bZ, cX+dZ),aX+bZ\Big),cX+dZ\right\rangle_g
\]
\[
=a\left\langle JN\Big(N(X, cX+dZ),aX+bZ\Big),cX+dZ\right\rangle_g
\]
\[
+b\left\langle JN\Big(N(Z, cX+dZ),aX+bZ\Big),cX+dZ\right\rangle_g
\]
\[
=ad\left\langle JN\Big(N(X, Z),aX+bZ\Big),cX+dZ\right\rangle_g
\]
\[
+bc\left\langle JN\Big(N(Z, X),aX+bZ\Big),cX+dZ\right\rangle_g
\]
\[
=ada\left\langle JN\Big(N(X, Z),X\Big),cX+dZ\right\rangle_g
\]
\[
+adb\left\langle JN\Big(N(X, Z),Z\Big),cX+dZ\right\rangle_g
\]
\[
+bca\left\langle JN\Big(N(Z, X),X\Big),cX+dZ\right\rangle_g
\]
\[
+bcb\left\langle JN\Big(N(Z, X),Z\Big),cX+dZ\right\rangle_g
\]
\[
=adad\left\langle JN\Big(N(X, Z),X\Big),Z\right\rangle_g
\]
\[
+adbc\left\langle JN\Big(N(X, Z),Z\Big),X\right\rangle_g
\]
\[
+bcad\left\langle JN\Big(N(Z, X),X\Big),Z\right\rangle_g
\]
\[
+bcbc\left\langle JN\Big(N(Z, X),Z\Big),X\right\rangle_g
\]
\[
=(a^2d^2-2abcd+b^2c^2)L(X,Z,X,Z)=(ad-bc)^2\ell(X,Z).
\]

Now let
\[
\ell_1(X,Z):=g(X,X)g(Z,Z)-g(X,Z)^2.
\]
Then
\[
\ell_1(X',Z')=g(aX+bZ,aX+bZ)g(cX+dZ.cX+dZ)-g(aX+bZ,cX+dZ)^2
\]
\[
=\big\{a^2g(X,X)+2abg(X,Z)+b^2g(Z,Z)\big\}
\big\{c^2g(X,X)+2cdg(X,Z)+d^2g(Z,Z)\big\}
\]
\[
-\big\{acg(X,X)+(ad+bc)g(X,Z)+bdg(Z,Z)\big\}^2
\]
\[
=a^2c^2g(X,X)^2+2a^2cdg(X,X)g(X,Z)+a^2d^2g(X,X)g(Z,Z)
\]
\[
+2abc^2g(X,Z)g(X,X)+4abcdg(X,Z)^2+2abd^2g(X,Z)g(Z,Z)
\]
\[
+b^2c^2g(Z,Z)g(X,X)+2b^2cdg(Z,Z)g(X,Z)+b^2d^2g(Z,Z)^2
\]
\[
-a^2c^2g(X,X)^2-(ad+bc)^2g(X,Z)^2-b^2d^2g(Z,Z)^2
\]
\[
-2ac(ad+bc)g(X,X)g(X,Z)-2(ad+bc)bdg(X,Z)g(Z,Z)
-2acbdg(X,X)g(Z,Z)
\]
\[
=+2a^2cdg(X,X)g(X,Z)+a^2d^2g(X,X)g(Z,Z)
\]
\[
+2abc^2g(X,Z)g(X,X)+4abcdg(X,Z)^2+2abd^2g(X,Z)g(Z,Z)
\]
\[
+b^2c^2g(Z,Z)g(X,X)+2b^2cdg(Z,Z)g(X,Z)
\]
\[
-(ad+bc)^2g(X,Z)^2
\]
\[
-2ac(ad+bc)g(X,X)g(X,Z)-2(ad+bc)bdg(X,Z)g(Z,Z)
\]
\[-2acbdg(X,X)g(Z,Z)
\]
\[
=(2a^2cd+2abc^2-2acad-2abc^2)g(X,X)g(X,Z)
(a^2d^2+b^2c^2-2abcd)g(X,X)g(Z,Z)
\]
\[
+(4abcd-a^2d^2-2abcd-b^2c^2)g(X,Z)^2
+(2abd^2+2b^2cd-2adbd-2bcbd)g(X,Z)g(Z,Z)
\]
\[
=(a^2d^2+b^2c^2-2abcd)g(X,X)g(Z,Z)
+(2abcd-a^2d^2-b^2c^2)g(X,Z)^2
\]
\[
=(ad-bc)^2\{g(X,X)g(Z,Z)-g(X,Z)^2\}
\]
\[=(ad-bc)^2\ell_1(X,Z).
\]

Therefore in general
\[
\ell(X',Z')=\left(\frac{\pa\{X',Z'\}}{\pa\{X,Z\}}\right)^2\ell(X,Z),
\]
\[
\ell_1(X',Z')=\left(\frac{\pa\{X',Z'\}}{\pa\{X,Z\}}\right)^2\ell_1(X,Z).
\]
and
\[
\frac{\ell(X',Z')}{\ell_1(X',Z')}=\frac{\ell(X,Z)}{\ell_1(X,Z)}.
\]
\qed

\vspace{0.1in}

$L(\cdot, \cdot\cdot,\cdots,\cdot\cdots)$ can be recaptured from 
$\ell(\cdot, \cdot\cdot)$, namely, we have the following theorem.
\begin{theorem}[Recapture $L$ from $\ell$]\label{2t-captures-4t-thm}
	\[
	L(X,Z,Y,W)=\frac14\frac{\pa^2}{\pa t\pa s}\Big|_{s=0,t=0}\ell(X+sY,Z+tW).
	\]
\end{theorem}

\proof Recall
\[
\ell(X,Z)=L(X,Z,X,Z)=\left\langle JN\Big(N(X,Z),X\Big),Z\right\rangle_g.
\]

In fact
\[
\frac{\pa^2}{\pa t\pa s}\Big|_{s=0,t=0}\ell(X+sY,Z+tW)
\]

\[
=\frac{\pa}{\pa t}\frac{\pa}{\pa s}
\left\langle JN\Big(N(X+sY,Z+tW),X+sY\Big),Z+tW\right\rangle_g\Big|_{s=0,t=0}
\]

\[
=\frac{\pa}{\pa t}
\left\langle  
JN\Big(N(Y,Z+tW),X+sY\Big),Z+tW
\right\rangle_g\Big|_{s=0,t=0}
\]
\[
+\frac{\pa}{\pa t}
\left\langle JN\Big(N(X+sY,Z+tW),Y\Big),Z+tW\right\rangle_g\Big|_{s=0,t=0}
\]

\[
=\frac{\pa}{\pa t}
\left\langle  
JN\Big(N(Y,Z+tW),X\Big),Z+tW
\right\rangle_g\Big|_{t=0}
\]
\[
+\frac{\pa}{\pa t}
\left\langle JN\Big(N(X,Z+tW),Y\Big),Z+tW\right\rangle_g\Big|_{t=0}
\]

\[
=
\left\langle  
JN\Big(N(Y,W),X\Big),Z+tW
\right\rangle_g\Big|_{t=0}
\]
\[
+\left\langle  
JN\Big(N(Y,Z+tW),X\Big),W
\right\rangle_g\Big|_{t=0}
\]
\[
+
\left\langle JN\Big(N(X,W),Y\Big),Z+tW\right\rangle_g\Big|_{t=0}
\]
\[
+
\left\langle JN\Big(N(X,Z+tW),Y\Big),W\right\rangle_g\Big|_{t=0}
\]

\[
=
\left\langle  
JN\Big(N(Y,W),X\Big),Z
\right\rangle_g
\]
\[
+
\left\langle  
JN\Big(N(Y,Z),X\Big),W
\right\rangle_g
\]
\[
+
\left\langle JN\Big(N(X,W),Y\Big),Z\right\rangle_g
\]
\[
+
\left\langle JN\Big(N(X,Z),Y\Big),W\right\rangle_g
\]

\[
=\left\langle JN\Big(N(X,Z),Y\Big),W\right\rangle_g
\]
\[
+
\left\langle  
JN\Big(N(Y,Z),X\Big),W
\right\rangle_g
\]
\[
+
\left\langle JN\Big(N(X,W),Y\Big),Z\right\rangle_g
\]
\[
+\left\langle  
JN\Big(N(Y,W),X\Big),Z
\right\rangle_g
\]
\[
=4L(X,Z,Y,W).
\]
\qed

\vspace{0.2in}

\section{Proof of the Other Main Result}\label{sec-proof-other-main results}

\proof[Proof of Theorem \ref{NN-is0-thm}]

If $\ell\equiv0$ then by Theorem (\ref{2t-captures-4t-thm}), $L\equiv0$. This implies
\[
JN\Big(N(X,Z),Y\Big)=0,
\]
and
\[
N^2(X,Z,Y)=N\Big(N(X,Z),Y\Big)=0
\]
for  smooth vector fields $X,Y,Z$, since $J$ is nonsingular.

\end{document}